%% file: cc.tex
\documentclass[12pt,twoside]{article}

\usepackage{amsgen,amsmath,amstext,amsbsy,amsopn,amsfonts,amssymb}

\usepackage{memoire}
\input{generique}
\def\j#1#2{\def\deux{#2} \ifx\deux\empty {\rm rk}\hskip .125em{{\goth #1}} \else {\rm rk}\hskip .125em{{\goth #1}_{#2}} \fi}
\def\k#1#2{\def\deux{#2} \ifx\deux\empty {\rm b}_{{\goth #1}} \else {\rm b}_{{\goth #1}_{#2}} \fi}

\def\an#1#2{\def\deux{#2} \ifx\deux\empty {\cal O}_{#1} \else {\cal O}_{#1,#2} \fi }
\def\han#1#2{\def\deux{#2} \ifx\deux\empty {\hat{{\cal O}}}_{#1} \else {\hat{{\cal O}}}_{#1,#2} \fi }

\def\r#1#2#3{\relax \def\un{#3} \relax \ifx\un\empty #1_{{\goth #2}} \else #1_{{\goth #2},{\goth #3}} \fi}

\def\dim{{\rm dim}\hskip .125em}

\def\dd{{\rm d}}
\def\rk{{\rm rk}\hskip .125em}

\def\r{{\rm r}}

\def\sy#1#2{{\rm S}^{#1}(#2)}

\def \ex #1#2{\mbox{\large$\wedge$}^{#1}(#2)}

\def\sl2{{\goth s}{\goth l}_{2}({\Bbb C})}

\renewcommand{\version}{{\rm Version r\'evis\'ee le~}\date} 
\renewcommand{\title}{\centerline{Complexe canonique
d'une alg\`ebre de Lie r\'eductive.}}
\renewcommand{\author}{Jean-Yves CHARBONNEL}
\renewcommand{\institution}{Universit\'e Paris 7 - CNRS}
\renewcommand{\lastname}{J-Y. Charbonnel}

\renewcommand{\entries}{%
\nextref{Bo} {N. Bourbaki} {Groupes et Alg\`ebres de Lie.
Chapitres 4, 5 et 6} {Masson. (1981), Paris, New-York,
Barcelone, Milan, Mexico, Rio de Janeiro}

\nextref{Bo1} {N. Bourbaki} {Groupes et Alg\`ebres de Lie.
Chapitres 7 et 8} {Diffusion CCLS. (1975), Paris}

\nextref{Ca}{R. W. Carter}{Finite Groups of Lie Type. Conjugacy
classes and complex characters} {Pure and Applied Mathematics. A
Wiley-Interscience Series of Texts, Monographs \& Tracts
(1985), Chichester, New-York, Brisbane, Toronto, Singapore}

\nextref{Dix}{J. Dixmier} {Champs de vecteurs adjoints sur les
groupes et alg\`ebres de Lie semi-simples} {Journal f\"ur die reine
und angewandte Mathematik {\bf Band. 309} (1979), p. 178--190}

\nextref{Ha}{R. Hartshorne} {Algebraic Geometry} {Graduate
Texts in Mathematics {\bf n$^{\circ}$ 52} (1977), Springer-Verlag, Berlin
Heidelberg New York}

\nextref{Mat}{H. Matsumura} {Commutative ring theory} {Cambridge
studies in advanced mathematics {\bf n$^{\circ}$ 8} (1986), Cambridge
University Press, Cambridge, London, New York, New Rochelle, Melbourne,
Sydney}

\nextref{Pa}{D. I. Panyushev} {The jacobian modules of a representation of a 
Lie algebra and geometry of commuting varieties} {Compositio Mathematica
{\bf n$^{\circ}$ 94} (1994), p. 181--199}

\nextref{Pr}{A. Premet} {Nilpotent commuting varieties of
reductive Lie algebras} {arXiv:math.RT/0302204 V1}

\nextref{Ric} {R. W. Richardson} {Commuting varieties of
semisimple Lie algebras and algebraic groups} {Compositio
Mathematica {\bf vol. 38} (1979), p. 311--322}

}

\renewcommand{\address}{Jean-Yves CHARBONNEL

Universit\'e Paris 7 - CNRS 

Institut de Math\'ematiques de Jussieu 

Th\'eorie des groupes

Case 7012

2 Place Jussieu

75251 Paris Cedex 05, France.

\makebox[1in][l]{jyc@math.jussieu.fr}

}

\begin{document}

\firstpage

\begin{abstract}
Let ${\goth g}$ be a finite dimensional complex reductive Lie algebra and
$\dv ..$ an invariant non degenerated bilinear form on 
${\goth g}\times {\goth g}$ which extends the Killing form of 
$[{\goth g},{\goth g}]$. We define the homology complex
$C_{\bullet}({\goth g})$. Its space is the algebra 
$\tk {{\Bbb C}}{\e Sg}\tk {{\Bbb C}}{\e Sg}\ex {}{{\goth g}}$ 
where $\e Sg$ and $\ex {}{{\goth g}}$ are the symmetric
and exterior algebras of ${\goth g}$. The differential of
$C_{\bullet}({\goth g})$ is the
$\tk {{\Bbb C}}{\e Sg}\e Sg$-derivation which associates to the
element $v$ of ${\goth g}$ the function 
$(x,y)\mapsto \dv v{[x,y]}$ on ${\goth g}\times {\goth g}$. Then the complex 
$C_{\bullet}({\goth g})$ has no homology in degree
strictly bigger than $\rk {\goth g}$.
\end{abstract}

\catcode`\@=11
\catcode`\;=\active\def;{\relax\ifhmode\ifdim\lastskip>\z@
\unskip\fi\kern.2em\fi\string;}
\catcode`\:=\active\def:{\relax\ifhmode\ifdim\lastskip>\z@\unsskip\fi
\penalty\@M\ \fi\string:}
\catcode`\!=\active\def!{\relax\ifhmode\ifdim\lastskip>\z@
\unskip\fi\kern.2em\fi\string!}
\catcode`\?=\active\def?{\relax\ifhmode\ifdim\lastskip>\z@
\unskip\fi\kern.2em\fi\string?}\frenchspacing\catcode`\@=12

\section{Introduction.} 
Soient ${\goth g}$ une alg\`ebre de Lie r\'eductive complexe de dimension
finie, $\e Sg$ et $\ex {}{{\goth g}}$ les alg\`ebres sym\'etrique et
ext\'erieure de ${\goth g}$. La repr\'esentation coadjointe de ${\goth g}$
s'identifie \`a sa repr\'esentation adjointe au moyen d'une forme
bilin\'eaire sym\'etrique invariante, non d\'eg\'en\'er\'ee, $\dv ..$, 
qui prolonge la forme de Killing de $[{\goth g},{\goth g}]$. On appelle 
{\it complexe canonique} de l'alg\`ebre de Lie ${\goth g}$ le complexe 
d'homologie $C_{\bullet}({\goth g})$ d'espace
$\tk {{\Bbb C}}{\e Sg}\tk {{\Bbb C}}{\e Sg}\ex {}{{\goth g}}$
dont la diff\'erentielle est la 
$\tk {{\Bbb C}}{\e Sg}\e Sg$-d\'erivation qui \`a
l'\'el\'ement $v$ de ${\goth g}$ associe la fonction 
$(x,y) \mapsto -\dv {v}{[x,y]}$. La vari\'et\'e
commutante de ${\goth g}$ est la sous-vari\'et\'e ${\goth C}_{{\goth g}}$ des
points $(x,y)$ de ${\goth g}\times {\goth g}$ tels que $[x,y]$ soit nul. Il
est facile de voir que le support dans ${\goth g}\times {\goth g}$ de 
l'homologie du complexe canonique de ${\goth g}$ est contenu dans la 
vari\'et\'e commutante de ${\goth g}$. Le r\'esultat principal de ce 
m\'emoire est le th\'eor\`eme:

\begin{Thd}
Le complexe $C_{\bullet}({\goth g})$ n'a pas d'homologie en
degr\'e strictement sup\'e\-rieur au rang  $\j g{}$ de ${\goth g}$. En outre, 
pour tout entier $i$ inf\'erieur \`a $\j g{}$, le complexe 
$C_{\bullet}({\goth g})$ a de l'homologie en degr\'e $i$.
\end{Thd}

On montrera dans un prochain article qu'un certain sous-complexe de 
$C_{\bullet}({\goth g})$ joue un r\^ole important dans l'\'etude de la
vari\'et\'e commutante de ${\goth g}$. Dans ce m\'emoire, le corps de base 
est le corps des nombres complexes, les espaces vectoriels et les alg\`ebres 
de Lie consid\'er\'es sont de dimension finie. On d\'esigne par ${\goth g}$ 
une alg\`ebre de Lie r\'eductive, par $G$ son groupe adjoint, par $\j g{}$ 
son rang et par $\dv ..$ une forme bilin\'eaire sym\'etrique, non 
d\'eg\'en\'er\'ee ${\goth g}$-invariante sur ${\goth g}$ qui prolonge la 
forme de Killing de ${\goth g}$. On identifie ${\goth g}$ \`a son dual au 
moyen de $\dv ..$. Si $V$ est un espace vectoriel, $V^{*}$ d\'esigne son dual.
Si ${\goth a}$ est une alg\`ebre de Lie et si ${\goth m}$ est un
sous-espace de ${\goth g}$, on note ${\goth m}(x)$ le sous-espace des
\'el\'ements de ${\goth m}$ qui centralisent l'\'el\'ement $x$ de
${\goth a}$ et ${\goth m}(x')$ le sous-espace des \'el\'ements de 
${\goth m}$ qui stabilisent la forme lin\'eaire $x'$ sur ${\goth a}$ pour
l'action coadjointe de ${\goth a}$ dans ${\goth a}^{*}$. On utilise la
topologie de Zariski sur les vari\'et\'es alg\'ebriques consid\'er\'ees. 
Si $X$ est une vari\'et\'e alg\'ebrique, $\an X{}$ d\'esigne le faisceau
structural de $X$ et $\an Xx$ l'anneau local au point $x$ de $X$. Sauf
mention contraire, un point de $X$ est un point ferm\'e. Conform\'ement \`a
l'usage, pour tout ouvert $Y$ de $X$ et pour tout faisceau ${\cal F}$ sur
$X$, $\Gamma (Y,{\cal F})$ d\'esigne l'espace des sections de ${\cal F}$
au dessus de $Y$. La suite de ce m\'emoire se divise en $3$ sections:
\begin{list}{}{}
\item 2) dimension projective et cohomologie, 
\item 3) exemple de complexe,
\item 4) complexe canonique d'une alg\`ebre de Lie r\'eductive.
\end{list}
\section{Dimension projective et cohomologie.} \label{p}
Soient $X$ une vari\'et\'e alg\'ebrique affine, de Cohen-Macaulay et 
${\Bbb C}[X]$ l'anneau des fonctions r\'eguli\`eres sur $X$. On note $U$ un
ouvert de $X$, $S$ le compl\'ementaire de $U$ dans $X$ et $p$ la
codimension de $S$ dans $X$. On d\'esigne par $P_{\bullet}$ un complexe
de ${\Bbb C}[X]$-modules projectifs de type fini, de longueur finie $l$, et
par $\varepsilon $ un morphisme d'augmentation du complexe
$P_{\bullet}$ d'image $K$, d'o\`u un complexe augment\'e de modules
$$ 0 \stackrel{\dd }\longrightarrow P_{l} \stackrel{\dd }\longrightarrow
P_{l-1} \stackrel{\dd }\longrightarrow \cdots
\stackrel{\dd }\longrightarrow P_{0} \stackrel{\varepsilon }
\longrightarrow K \rightarrow 0 \mbox{ .}$$

On note ${\cal P}_{\bullet}$ la localisation sur $X$ du complexe
$P_{\bullet}$, ${\cal K}_{0}$ la localisation sur $X$ du noyau de
$\varepsilon $, ${\cal K}$ la localisation sur $X$ de $K$ et
${\cal K}_{i}$ le noyau du morphisme 
${\cal P}_{i}\rightarrow {\cal P}_{i-1}$ pour $i$ entier strictement
positif.

\begin{lemme} \label{lp}
On suppose que $S$ contient le support de l'homologie du complexe
augment\'e $P_{\bullet}$.

{\rm i)} Pour tout $\an X{}$-module projectif ${\cal P}$, le groupe de
cohomologie ${\rm H}^{i}(U,{\cal P})$ est nul pour tout entier $i$
strictement positif et strictement inf\'erieur \`a $p-1$.

{\rm ii)} Pour tout entier naturel $j$ inf\'erieur \`a $l$, le groupe de
cohomologie  ${\rm H}^{i}(U,{\cal K}_{l-j})$ est nul si $i$ est un entier
strictement positif et strictement inf\'erieur \`a $p-j$.

{\rm iii)} Le groupe de cohomologie ${\rm H}^{i}(U,{\cal K})$ est nul pour
tout entier $i$ strictement positif et strictement inf\'erieur \`a $p-l-1$.
\end{lemme}

\begin{proof}
i) Soit $i$ un entier positif strictement inf\'erieur \`a $p-1$. Puisque le
foncteur  H$^{i}(U,\bullet )$ commute \`a la somme directe, il suffit de
montrer que H$^{i}(U,\an X{})$ est nul. Puisque $S$ est un ferm\'e de $X$,
on a la suite exacte longue de cohomologie relative
$$ \cdots \rightarrow 
{\rm H}^{i}_{S}(X,\an X{}) \rightarrow {\rm H}^{i}(X,\an X{})
\rightarrow {\rm H}^{i}(U,\an X{}) \rightarrow 
{\rm H}^{i+1}_{S}(X,\an X{}) \rightarrow \cdots \mbox{ ;}$$
or H$^{i}(X,\an X{})$ est nul car $X$ est affine; donc  H$^{i}(U,\an X{})$ est
isomorphe au groupe de cohomologie relative H$_{S}^{i+1}(X,\an X{})$ pour
tout entier strictement positif $i$. La vari\'et\'e $X$ \'etant de
Cohen-Macaulay, la codimension $p$ de $S$ est \'egale \`a la profondeur
de son id\'eal de d\'efinition dans ${\Bbb C}[X]$; or $i+1$ est strictement
inf\'erieur \`a $p$; donc H$_{S}^{i+1}(X,\an X{})$ est nul. Il en r\'esulte 
que H$^{i}(U,\an X{})$ est nul. 

ii) Puisque $S$ contient le support de l'homologie du complexe
$P_{\bullet}$, pour tout entier naturel $j$, on a la suite exacte courte de
$\an U{}$-modules
$$ 0 \rightarrow 
{\cal K}_{j+1}\left \vert \right. _{U} \rightarrow {\cal P}_{j+1}\left
\vert \right. _{U} \rightarrow {\cal K}_{j} \left \vert \right. _{U}
\rightarrow 0 \mbox{ .}$$
Il en r\'esulte la suite exacte longue de cohomologie
$$ \cdots \rightarrow {\rm H}^{i}(U,{\cal P}_{j+1})
\rightarrow {\rm H}^{i}(U,{\cal K}_{j}) \rightarrow 
{\rm H}^{i+1}(U,{\cal K}_{j+1}) \rightarrow 
{\rm H}^{i+1}(U,{\cal P}_{j+1}) \rightarrow \cdots \mbox{ ;}$$
donc d'apr\`es (i), pour $i$ entier strictement positif et strictement
inf\'erieur \`a $p-2$, les groupes de cohomologie
H$^{i}(U,{\cal K}_{j})$ et H$^{i+1}(U,{\cal K}_{j+1})$ sont isomorphes car
$P_{j+1}$ est un module projectif. Puisque ${\cal P}_{i}$
est nul pour $i$ stictement sup\'erieur \`a $l$, ${\cal K}_{l-1}$ et 
${\cal P}_{l}$ ont des restrictions \`a $U$ isomorphes. En particulier, 
H$^{i}(U,{\cal K}_{l-1})$ est nul pour $i$ strictement inf\'erieur \`a
$p-1$ d'apr\`es (i). Un raisonnement par r\'ecurrence sur $j$ montre alors
que le groupe de cohomologie H$^{i}(U,{\cal K}_{l-j})$ est nul pour $i$
strictement inf\'erieur \`a $p-j$.

iii) Soit $i$ un entier strictement positif et strictement inf\'erieur \`a
$p-l-1$. De la suite exacte courte
$$ 0 \rightarrow {\cal K}_{0}\left \vert \right. _{U}
\rightarrow {\cal P}_{0}\left \vert \right. _{U} \rightarrow 
{\cal K} \left \vert \right. _{U} \rightarrow 0 \mbox{ ,}$$
on tire la suite exacte
$$ {\rm H}^{i}(U,{\cal P}_{0}) \rightarrow {\rm H}^{i}(U,{\cal K})
\rightarrow {\rm H}^{i+1}(U,{\cal K}_{0}) \rightarrow 
{\rm H}^{i+1}(U,{\cal P}_{0}) \mbox{ ;}$$
or d'apr\`es (i), les groupes H$^{i}(U,{\cal P}_{0})$ et 
H$^{i+1}(U,{\cal P}_{0})$ sont nuls; donc d'apr\`es (ii), H$^{i}(U,{\cal K})$
est nul. 
\end{proof}

On suppose que $K$ est le noyau d'un morphisme $\theta $ du module
projectif de type fini $P$ dans le module projectif de type fini $P'$.
Soient $R$ l'image de $\theta $ et $R'$ un sous-module de $P'$ qui
contient $R$. On d\'esigne respectivement par ${\cal P}$, ${\cal P}'$,
${\cal R}$, ${\cal R}'$, $\theta $ les localisations respectives sur $X$ de
$P$, $P'$, $R$, $R'$, $\theta $. 

\begin{prop} \label{pp}
On suppose que $p$ est strictement sup\'erieur \`a $l+2$, que $X$ est
normal et que $S$ contient le support de l'homologie du complexe
augment\'e $P_{\bullet}$.

{\rm i)} Si ${\cal R}$ et ${\cal R}'$ ont m\^eme restriction \`a $U$, alors
$R$ est \'egal \`a $R'$.  

{\rm ii)} Le complexe $P_{\bullet}$ est une r\'esolution de $K$. 

{\rm iii)} Si ${\cal R}$ et ${\cal R}'$ ont m\^eme restriction \`a $U$, le
complexe augment\'e
$$ 0 \rightarrow P_{l} \rightarrow \cdots \rightarrow P_{0}
\rightarrow P \rightarrow R' \rightarrow 0 \mbox{ ,}$$
est une r\'esolution projective de $R'$ de longueur $l+1$.
\end{prop}

\begin{proof}
i) On suppose que ${\cal R}$ et ${\cal R}'$ ont m\^eme restriction \`a $U$.
Alors on a la suite exacte courte de $\an U{}$-modules
$$ 0 \rightarrow {\cal K}\left \vert \right. _{U} \rightarrow 
{\cal P}\left \vert \right. _{U} \rightarrow 
{\cal R}'\left \vert \right. _{U}\rightarrow 0 \mbox{ ,}$$ 
d'o\`u la suite exacte longue de cohomologie
$$ 0 \rightarrow \Gamma (U,{\cal K}) \rightarrow \Gamma (U,{\cal P})
\rightarrow \Gamma (U,{\cal R}') \rightarrow {\rm H}^{1}(U,{\cal K})
\rightarrow \cdots \mbox{ .}$$ 
D'apr\`es l'assertion (iii) du lemme \ref{lp}, H$^{1}(U,{\cal K})$ est nul
car $1$ est strictement inf\'erieur \`a $p-l-1$; donc on a la suite exacte
courte 
$$ 0 \rightarrow \Gamma (U,{\cal K}) \rightarrow \Gamma (U,{\cal P})
\rightarrow \Gamma (U,{\cal R}') \rightarrow 0 \mbox{ .}$$ 
Puisque le compl\'ementaire de $U$ dans $X$ est de codimension
sup\'erieure \`a $2$, $U$ est partout dense dans $X$ car $X$ est une 
vari\'et\'e de Cohen-Macaulay. En outre, la restriction de $X$ \`a $U$ est un 
isomorphisme de $P$ sur $\Gamma (U,{\cal P})$ car $X$ est une vari\'et\'e
normale. Soit $\varphi $ dans $R'$. Alors il existe un \'el\'ement $\psi $
de $P$ dont l'image par $\theta $ a m\^eme restriction \`a
$U$ que $\varphi $. Il en r\'esulte que $\varphi -\theta (\psi )$ est un
\'el\'ement de torsion de $P'$; donc $\varphi $ est \'egal \`a
$\theta (\psi )$. Par suite, $R$ est \'egal \`a $R'$ car $R'$ contient
l'image $R$ de $\theta $. 

ii) Pour $j$ entier naturel, on note $Z_{j}$ l'espace des cycles
de degr\'e $j$ du complexe $P_{\bullet}$ et pour $j=-1$, on pose
$K=Z_{-1}$. Par hypoth\`ese, pour tout entier $j$ sup\'erieur \`a $-1$,
$S$ contient le support de l'homologie du complexe augment\'e
$$ 0 \rightarrow P_{l} \rightarrow \cdots \rightarrow P_{j+1}
\rightarrow Z_{j} \rightarrow 0 \mbox{ ;}$$
or pour tout $j$ sup\'erieur \`a $-1$, $p$ est strictement sup\'erieur
\`a $l-j$; donc d'apr\`es l'assertion (i), $Z_{j}$ est l'image de $P_{j+1}$
dans $P_{j}$ pour tout entier naturel $j$ et $Z_{-1}$ est l'image de
$P_{0}$. Cela revient \`a dire que le complexe $P_{\bullet}$ est une
r\'esolution de $K$.

iii) L'assertion r\'esulte des assertions (i) et (ii).
\end{proof}

\begin{cor} \label{cp}
Soient ${\cal C}_{\bullet}$ un complexe de  ${\Bbb C}[X]$-modules de type
fini de longueur strictement positive $l$ et $d$ un entier naturel. On suppose
que les trois conditions suivantes sont r\'ealis\'ees:
\begin{list}{}{}
\item {\rm 1)} $X$ est normal et $S$ contient le support de l'homologie du
complexe $C_{\bullet}$,
\item {\rm 2)} pour tout entier naturel $i$, $C_{i}$ est un sous-module
d'un module libre,
\item {\rm 3)} pour tout entier strictement positif $i$, la dimension
projective de $C_{i}$ est inf\'erieure \`a $i+d$.
\end{list}
Pour $j$ entier naturel inf\'erieur \`a $l$, on d\'esigne par $Z_{j}$
l'espace des cycles de degr\'e $j$ du complexe $C_{\bullet}$. Si $p$ est
strictement sup\'erieur \`a $2l+d$, alors le complexe $C_{\bullet}$
est acyclique. En outre, la dimension projective de $Z_{j}$ est inf\'erieure
\`a $2l+d-j-1$.
\end{cor}

\begin{proof}
On suppose $p$ strictement sup\'erieur \`a $2l+d$. On
montre en raisonnant par r\'ecurrence sur $l-j$ que le complexe
$$ 0 \rightarrow C_{l} \rightarrow \cdots \rightarrow C_{j+1}
\rightarrow Z_{j} \rightarrow 0 \mbox{ ,}$$
est acyclique et que la dimension projective de $Z_{j}$ est inf\'erieure
\`a $2l+d-j-1$. Pour $j=l$, le support de $Z_{j}$ est contenu dans $S$;
or $C_{l}$ est un sous-module d'un module libre; donc $Z_{l}$ est nul.
Puisque $Z_{l-1}$ est un sous-module d'un module libre,
$Z_{l-1}$ est l'espace des bords de degr\'e $l-1$ du complexe
$C_{\bullet}$ d'apr\`es l'assertion (iii) de la proposition \ref{pp} car 
$l+d+1$ est strictement inf\'erieur \`a $p$. Il en r\'esulte que la dimension
projective de $Z_{l-1}$ est inf\'erieure \`a $l+d$. On suppose
$j$ inf\'erieur \`a $l-2$ et l'assertion vraie pour $j+1$. Par
hypoth\`ese, $C_{j+1}$ a une r\'esolution projective $P_{\bullet}$ de
longueur au plus $j+d+1$. D'apr\`es l'hypoth\`ese de r\'ecurrence,
$Z_{j+1}$ a une r\'esolution projective $Q_{\bullet}$ de longueur
au plus $2l+d-j-2$. On en d\'eduit un complexe augment\'e 
$R_{\bullet}$ de modules projectifs de longueur inf\'erieure \`a
$2l+d-j-1$,
$$ 0 \rightarrow Q_{2l+d-j-2} \oplus P_{2l+d-j-1} \rightarrow 
\cdots  \rightarrow Q_{0} \oplus P_{1}
\rightarrow P_{0} \rightarrow  Z_{j} \rightarrow 0 \mbox{ .}$$
D\'esignant par $\dd $ les diff\'erentielles des complexes $Q_{\bullet}$
et $P_{\bullet}$, la restriction \`a $Q_{i}\oplus P_{i+1}$ de la
diff\'erentielle de $R_{\bullet}$ est l'application
$$(x,y) \mapsto (\dd x,\dd y +(-1)^{i} \delta (x)) \mbox{ ,}$$ 
o\`u $\delta $ est l'application de $Q_{i}$ dans $P_{i}$ qui r\'esulte de
l'injection de $Z_{j+1}$ dans $C_{j+1}$. Puisque $P_{\bullet}$ et
$Q_{\bullet}$ sont des r\'esolutions projectives, le complexe
$R_{\bullet}$ est un complexe de modules projectifs qui n'a pas
d'homologie en degr\'e positif; donc le support de l'homologie du complexe 
augment\'e $R_{\bullet}$ est contenu dans $S$. Il r\'esulte alors de 
l'assertion (iii) de la proposition \ref{pp} que $R_{\bullet}$ est une 
r\'esolution projective de longueur $2l+d-j-1$ de $Z_{j}$ car $Z_{j}$ est un 
sous-module d'un module libre, d'o\`u le corollaire.
\end{proof}
\section{Exemple de complexe.} \label{co}
Dans cette section, on donne un exemple de complexe d'homologie
qui sera utilis\'e dans la suite de ce m\'emoire. Soient $X$
une vari\'et\'e alg\'ebrique affine, irr\'eductible et ${\Bbb C}[X]$
l'alg\`ebre des fonctions r\'eguli\`eres sur $X$.  Dans ce qui suit, on 
d\'esigne par $V$ un espace vectoriel. On note respectivement $\es SV$ et 
$\ex {}V$ les alg\`ebres sym\'etrique et ext\'erieure de $V$. Pour tout 
entier $i$, $\sy iV$ et $\ex iV$ d\'esignent respectivement les sous-espaces 
de degr\'e $i$ pour leurs graduations usuelles des alg\`ebres $\es SV$ et 
$\ex {}V$. En particulier, pour $i$ strictement n\'egatif, $\sy iV$ et 
$\ex iV$ sont nuls. 

\subsection{} Soient $\lambda $ une application lin\'eaire de $V$ dans
${\Bbb C}[X]$ et $\theta _{\lambda }$ l'application lin\'eaire 
$$\tk {{\Bbb C}}{{\Bbb C}[X]}V \stackrel{\theta _{\lambda }}
\longrightarrow {\Bbb C}[X] \mbox{ , } a\tens v \mapsto a\lambda (v) 
\mbox{ ,}$$ 
o\`u $a$ est dans ${\Bbb C}[X]$ et o\`u $v$ est dans $V$.

\begin{Def} \label{dco1}
On appelle complexe canonique associ\'e \`a $\lambda $ et on le note
$C_{\bullet}(\lambda )$ le complexe d'espace
$\tk {{\Bbb C}}{{\Bbb C}[X]}\ex {}V$ dont la diff\'erentielle est la
${\Bbb C}[X]$-d\'erivation de l'alg\`ebre 
$\tk {{\Bbb C}}{{\Bbb C}[X]}\ex {}V$ qui prolonge $\theta _{\lambda }$. La
graduation naturelle de $\ex {}V$ induit sur $C_{\bullet}(\lambda )$ une
structure de complexe  gradu\'e. On d\'esigne par $K_{\lambda }$ et 
$I_{\lambda }$ le noyau et l'image de $\theta _{\lambda }$. 
\end{Def}

Pour $\pi $ automorphisme de la vari\'et\'e $X$, on d\'esigne par
$\pi ^{\#}$ l'automorphisme de l'alg\`ebre 
$\tk {{\Bbb C}}{{\Bbb C}[X]}\ex {}V$ qui \`a l'application $\varphi $ de $X$
dans $\ex {}V$ associe $\varphi \rond \pi $ et $\lambda _{\pi }$ l'application
de $V$ dans ${\Bbb C}[X]$ qui \`a l'\'el\'ement $v$ associe la fonction
$\lambda (v)\rond \pi $ sur $X$.  

\begin{lemme}\label{lco1}
Soient $\pi $ un automorphisme de $X$ et $\dd $ la diff\'erentielle du
complexe $C_{\bullet}(\lambda )$.

{\rm i)} Le module $K_{\lambda _{\pi }}$ et l'id\'eal 
$I_{\lambda _{\pi }}$ sont les images respectives de $K_{\lambda }$ et de
$I_{\lambda }$ par $\pi ^{\#}$.

{\rm ii)} La diff\'erentielle du complexe $C_{\bullet}(\lambda _{\pi })$
est l'application $\pi ^{\#}\rond \dd \rond (\pi ^{\#})^{-1}$.

{\rm iii)} L'homologie du complexe $C_{\bullet}(\lambda )$ est un
\hbox{${\Bbb C}[X]$-module} gradu\'e de type fini. En outre, pour 
tout entier naturel $j$, le support du $j$-i\`eme groupe d'homologie de ce
complexe est une partie ferm\'ee  de $X$.

{\rm iv)} Pour tout entier naturel $j$, l'image par $\pi ^{-1}$ du
support du \hbox{$j$-i\`eme} groupe d'homologie du complexe
$C_{\bullet}(\lambda )$ est le support du $j$-i\`eme groupe d'homologie du
complexe $C_{\bullet}(\lambda _{\pi })$.
\end{lemme}

\begin{proof}
i) Soit $\varphi $ dans $\tk {{\Bbb C}}{{\Bbb C}[X]}V$. L'\'el\'ement
$\pi ^{\#}(\varphi )$ appartient \`a $K_{\lambda _{\pi }}$ si et seulement
si $\lambda (\varphi \rond \pi (x))(\pi (x))$ est nul pour tout $x$ dans $X$.
Puisque $\pi $ est un automorphisme, cela revient \`a dire que 
$K_{\lambda }$ contient $\varphi $; donc $K_{\lambda _{\pi }}$ est l'image de
$K_{\lambda }$ par $\pi ^{\#}$. L'image de $\varphi $ par le
morphisme $\theta _{\lambda _{\pi }}$ est la fonction 
$x\mapsto \lambda (\varphi (x))(\pi (x))$ sur $X$. C'est en particulier l'image
de $\varphi \rond \pi ^{-1}$ par $\pi ^{\#}\rond \theta _{\lambda }$; donc
$I_{\lambda _{\pi }}$ est \'egal \`a $\pi ^{\#}(I_{\lambda })$. 

ii) Pour tout $v$ dans $V$ et pour tout $\varphi $ dans ${\Bbb C}[X]$, on a
$$ \pi ^{\#}\rond \dd \rond (\pi ^{\#})^{-1}(\varphi \tens v) = 
\varphi \lambda _{\pi }(v)  \mbox{ ;}$$
or l'application $\pi ^{\#}\rond \dd \rond (\pi ^{\#})^{-1}$ est une
d\'erivation ${\Bbb C}[X]$-lin\'eaire de l'alg\`ebre \sloppy 
\hbox{$\tk {{\Bbb C}}{{\Bbb C}[X]}\ex {}V$}; donc d'apr\`es l'unicit\'e du
prolongement des d\'erivations, $\pi ^{\#}\rond \dd \rond (\pi ^{\#})^{-1}$
est la diff\'erentielle du complexe $C_{\bullet}(\lambda _{\pi })$.

iii) Puisque $C_{j}(\lambda )$ est nul pour $j$ strictement
sup\'erieur \`a $\dim V$, $C_{\bullet}(\lambda )$ est un 
${\Bbb C}[X]$-module de type fini car $C_{j}(\lambda )$ est un 
${\Bbb C}[X]$-module de type fini pour tout entier $j$. La diff\'erentielle
$\dd $ \'etant ${\Bbb C}[X]$-lin\'eaire, l'homologie du complexe 
$C_{\bullet}(\lambda )$ est un sous-quotient du  
${\Bbb C}[X]$-module $C_{\bullet}(\lambda )$; donc l'homologie du
complexe $C_{\bullet}(\lambda )$ est un module de type fini. Par suite, 
pour tout entier naturel $j$, le $j$-i\`eme groupe d'homologie  
H$_{j}(C_{\bullet}(\lambda ))$ du complexe $C_{\bullet}(\lambda )$ est un
\hbox{${\Bbb C}[X]$-module} de type fini. Il en r\'esulte que le support
du module H$_{j}(C_{\bullet}(\lambda ))$ est ferm\'e dans $X$.

iv) Soient $j$ un entier naturel, $J_{j}$ et $J_{j,\pi }$ les id\'eaux de
d\'efinition dans ${\Bbb C}[X]$ des supports des $j$-i\`emes groupes
d'homologie des complexes 
$C_{\bullet}(\lambda )$ et $C_{\bullet}(\lambda _{\pi })$. Si $a$
est un cycle de degr\'e $j$ de $C_{\bullet}(\lambda )$ et si $p$ est
dans $J_{j}$, pour $m$ entier naturel assez grand, $p^{m}a$ est un bord de
$C_{\bullet}(\lambda )$; donc d'apr\`es l'assertion (ii),
$\pi ^{\#}(p)^{m}\pi ^{\#}(a)$ est un bord de 
$C_{\bullet}(\lambda _{\pi })$. Il r\'esulte alors de cette assertion
que $J_{j,\pi }$ contient $\pi ^{\#}(J_{j})$. De m\^eme,
$J_{j}$ contient $(\pi ^{\#})^{-1}(J_{j,\pi })$; donc $J_{j,\pi }$ est
l'image par $\pi ^{\#}$ de $J_{j}$. Par suite, le support du $j$-i\`eme
groupe d'homologie de $C_{\bullet}(\lambda _{\pi })$ est l'image par 
$\pi ^{-1}$ du support du $j$-i\`eme groupe d'homologie de
$C_{\bullet}(\lambda )$.
\end{proof}

\subsection{} Soient
$W$ un espace vectoriel et $\tau $ une application r\'eguli\`ere de $X$
dans l'espace des applications lin\'eaires de $V$ dans $W$. 

\begin{Def}\label{dco2}
On appelle complexe canonique associ\'e \`a $\tau $ et on le note 
$C_{\bullet}(\tau )$ le complexe $C_{\bullet}(\lambda )$, d\'efini en
\ref{dco1}, o\`u $\lambda $ est l'application lin\'eaire de $V$ dans 
${\Bbb C}[X\times W^{*}]$ qui \`a $v$ associe la fonction
$(x,w') \mapsto \dv {w'}{\tau (x)(v)}$. 
\end{Def}

On rappelle que $W^{*}$ est le dual de $W$ et que l'alg\`ebre 
${\Bbb C}[X\times W^{*}]$ est canoniquement isomorphe \`a l'alg\`ebre 
$\tk {{\Bbb C}}{{\Bbb C}[X]}\es SW$. 

\begin{lemme}\label{lco2}
On note ${\goth C}_{\tau }$ la vari\'et\'e des z\'eros de $I_{\lambda }$.

{\rm i)} La vari\'et\'e ${\goth C}_{\tau }$ est l'ensemble des \'el\'ements
$(x,w')$ de $X\times W^{*}$ qui satisfont la condition  suivante: $w'$ est
orthogonal \`a l'image de $\tau (x)$. 

{\rm ii)} Le radical de $I_{\lambda }$ est l'ensemble des \'el\'ements
$\varphi $ de $\tk {{\Bbb C}}{{\Bbb C}[X]}\es SW$ qui satisfont la
condition suivante: pour tout $x$ dans $X$, $\varphi (x)$ appartient \`a
l'id\'eal de $\es SW$ engendr\'e par l'image de $\tau (x)$. 

{\rm iii)} Le support de l'homologie du complexe $C_{\bullet}(\tau )$ est
contenu dans ${\goth C}_{\tau }$.

{\rm iv)} Tout bord du complexe $C_{\bullet}(\tau )$ est nul sur 
${\goth C}_{\tau }$.
\end{lemme}

\begin{proof}
i) Par d\'efinition, ${\goth C}_{\tau }$ contient $(x,w')$ si et seulement si
\sloppy $\dv {w'}{\tau (x)(v)}$ est nul pour tout $v$ dans $V$, d'o\`u
l'assertion. 

ii) Pour tout $x$ dans $W$, l'id\'eal de $\es SW$ engendr\'e par l'image de
$\tau (x)$ est l'id\'eal de d\'efiniton de l'orthogonal dans $W^{*}$ de
l'image de  $\tau (x)$; donc $\varphi $ appartient au radical de
$I_{\lambda }$ si et seulement si pour tout $x$ dans $X$, $\varphi (x)$
appartient \`a l'id\'eal de $\es SW$ engendr\'e par l'image de $\tau (x)$.

iii) On note ${\cal C}_{\bullet}(\tau )$ la localisation sur 
$X\times W^{*}$ du complexe $C_{\bullet}(\tau )$. Soit
$(x,w')$ un point de $X\times W^{*}$ qui n'est pas dans 
${\goth C}_{\tau }$. Puisque la forme lin\'eaire $w'\rond \tau (x)$ n'est
pas nulle, il existe un ouvert affine $U$ de $X\times W^{*}$ qui
contient $(x,w')$ et des applications r\'eguli\`eres 
$\poi {\varphi }1{,\ldots,}{n}{}{}{}$ de $U$ dans $V$ qui satisfont
les conditions suivantes: pour tout $(z,z')$ dans $U$, 
$\dv {z'}{\tau (z)(\varphi _{n}(z,z'))}$ est \'egal
\`a $1$ et $\poi {z,z'}{}{,\ldots,}{}{\varphi }{1}{n-1}$ est une base du
noyau de $z'\rond \tau (z)$. On d\'esigne par $C_{\bullet}$
l'espaces des sections au dessus de $U$ du complexe
${\cal C}_{\bullet}(\tau )$. Puisque le bord de $\varphi_{n}$ est \'egal \`a 
$1$, $c$ est le bord de $\varphi _{n}\wedge c$ pour tout cycle $c$ de 
$C_{\bullet}$; donc $C_{\bullet}$ est acyclique. Vu l'arbitraire de $(x,w')$ 
et de l'ouvert affine $U$ contenant $(x,w')$, ${\goth C}_{\tau }$ contient le 
support de l'homologie du complexe $C_{\bullet}(\tau )$.

iv) L'espace des bords de $C_{\bullet}(\tau )$ est contenu dans l'id\'eal de
\sloppy  \hbox{$\tk {{\Bbb C}}{{\Bbb C}[X]}\tk {{\Bbb C}}{\es SW}\ex {}V$} 
engendr\'e par $\lambda (V)$, d'o\`u l'assertion car $\lambda (V)$ engendre 
l'id\'eal $I_{\lambda }$ par d\'efinition.
\end{proof}

\section{Complexe canonique d'une alg\`ebre de Lie r\'eductive.}
\label{ca} 
Soit $\lambda _{{\goth g}}$ l'application lin\'eaire de ${\goth g}$ dans 
$\tk {{\Bbb C}}{{\goth g}}{\goth g}$ qui \`a l'\'el\'ement $v$ associe la 
fonction $(x,y) \mapsto \dv {v}{[x,y]}$. On utilise alors les notations de la 
d\'efinition \ref{dco1}. On d\'esigne par $I_{{\goth g}}$ l'id\'eal de 
$\tk {{\Bbb C}}{\e Sg}\e Sg$ engendr\'e par $\lambda _{{\goth g}}({\goth g})$.

\begin{Def} \label{dca}
On appelle complexe canonique de l'alg\`ebre de Lie ${\goth g}$
le complexe canonique associ\'e \`a $\lambda _{{\goth g}}$ au sens de
\ref{dco1}.
\end{Def}

\subsection{} On note ${\goth C}_{{\goth g}}$ la vari\'et\'e commutante
de ${\goth g}$. Par d\'efinition, ${\goth C}_{{\goth g}}$ est la vari\'et\'e
des z\'eros dans ${\goth g}\times {\goth g}$ de l'id\'eal $I_{{\goth g}}$.

\begin{lemme}\label{lca2}
{\rm i)} Le support de l'homologie du complexe $C_{\bullet}({\goth g})$ est 
contenu dans ${\goth C}_{{\goth g}}$.

{\rm ii)} Tout bord du complexe $C_{\bullet}({\goth g})$ est nul sur 
${\goth C}_{{\goth g}}$. 
\end{lemme}

\begin{proof}
Le lemme r\'esulte des assertions (iii) et (iv) du lemme \ref{lco2}.
\end{proof}
\subsection{} 
On rappelle que d'apr\`es \cite{Ric}, la vari\'et\'e commutante 
${\goth C}_{{\goth g}}$ de ${\goth g}$ est irr\'eductible et de dimension 
$\dim {\goth g}+\j g{}$.

\begin{Th} \label{tca}
Le complexe canonique $C_{\bullet}({\goth g})$ de ${\goth g}$ n'a pas
d'homologie en degr\'e strictement sup\'erieur \`a $\rk {\goth g}$. En outre,
pour tout entier $i$, inf\'erieur \`a $\j g{}$, le complexe 
$C_{\bullet}({\goth g})$ a de l'homologie en degr\'e $i$. 
\end{Th}

\begin{proof}
Pour $i$ entier strictement positif, on note $Z_{i}({\goth g})$ l'espace des 
cycles de degr\'e $i$ du complexe $C_{\bullet}({\goth g})$. Alors pour $i$
entier strictement positif et strictement inf\'erieur \`a ${\goth g}$, on a un 
complexe augment\'e
$$ 0 \rightarrow C_{\dim {\goth g}}({\goth g}) \rightarrow \cdots 
\rightarrow C_{i+1}({\goth g})\rightarrow Z_{i}({\goth g})\rightarrow 0
\mbox{ .}$$
D'apr\`es l'assertion (i) du lemme \ref{lca2}, le support de l'homologie de
ce complexe est contenu dans ${\goth C}_{{\goth g}}$; or la codimension de 
${\goth C}_{{\goth g}}$ dans ${\goth g}\times {\goth g}$ est \'egale \`a 
$\dim {\goth g}-\j g{}$ et pour tout entier $i$ strictement positif, 
$C_{i}({\goth g})$ est un module libre; donc d'apr\`es l'assertion (iii) de la
proposition \ref{pp}, pour $i$ strictement sup\'erieur \`a $\j g{}$ et 
strictement inf\'erieur \`a $\dim {\goth g}$, le complexe ci-dessus est 
acyclique. Il en r\'esulte que le complexe $C_{\bullet}({\goth g})$ n'a pas
d'homologie en degr\'e strictement sup\'erieur \`a $\j g{}$. 

Soit $L_{{\goth g}}$ le sous-module des \'el\'ements $\varphi $ de 
$\tk {{\Bbb C}}{\e Sg}{\goth g}$ qui satisfont l'\'egalit\'e
$[x,\varphi (x)]=0$ pour tout $x$ dans ${\goth g}$. D'apr\`es 
\cite{Dix}(\S 2), $L_{{\goth g}}$ est un module libre de rang $\j g{}$; donc
pour tout entier naturel $i$, inf\'erieur \`a $\j g{}$, le module
$\ex i{L_{{\goth g}}}$ n'est pas nul. Soient $i$ un entier naturel, 
inf\'erieur \`a $\j g{}$ et $\psi $ un \'el\'ement non nul de 
$\ex {i}{L_{{\goth g}}}$. Pour tout $\varphi $ dans $L_{{\goth g}}$, 
l'application $(x,y) \mapsto \varphi (x)$ est un cycle du 
complexe $C_{\bullet}({\goth g})$; donc l'application 
$(x,y)\mapsto \psi (x)$ est cycle de degr\'e $i$ du complexe 
$C_{\bullet}({\goth g})$. En outre, ce cycle ne s'annule pas sur 
${\goth C}_{{\goth g}}$ car $\psi $ n'est pas nul; donc d'apr\`es 
l'assertion (ii) du lemme \ref{lca2}, $C_{\bullet}({\goth g})$ a de 
l'homologie en degr\'e $i$.
\end{proof}

\references

\lastpage

\end{document}

%% file: generique
\def\poi#1#2#3#4#5#6#7{\def\un{#5#6#7}\def\deux{#6#7}
\def\trois{#2#4} \def\cinq{#3#4#5}
\ifx\un\empty {#1}_{#2}{#3}{#1}_{#4} \else
\ifx\deux\empty {#5}(#1_{#2}){#3}{#5}(#1_{#4}) \else
\ifx\trois\empty {#5}_{#6}(#1){#3}{#5}_{#7}(#1) \else
{#5_{#6}}(#1_{#2}){#3}{#5_{#7}}(#1_{#4}) \fi \fi \fi}
\def\rond{\raisebox{.3mm}{\scriptsize$\circ$}}

\def\tens{\raisebox{.3mm}{\scriptsize$\otimes$}}

\def\dv#1#2{\langle {#1},{#2}\rangle}
    
\def\tk#1#2{{#2}\otimes _{#1}}

\def\ec#1#2#3#4#5{\def\un{#3#4#5}\def\deux{#3#5}\def\trois{#3}
\def\four{#2#4#5}\def\five{#2#5}\def\six{#2}\def\seven{#3#4}
\def\eight{#2#4} \def\nine{#2#3#4}
\ifx\nine\empty {\rm #1}_{#5} \else
\ifx\un\empty {\rm #1}({\goth #2}) \else
\ifx\deux\empty {\rm #1}({\goth #2}_{#4}) \else
\ifx\trois\empty {\rm #1}_{#5}({\goth #2}_{#4}) \else
\ifx\four\empty {\rm #1}(#3) \else
\ifx\five\empty {\rm #1}(#3_{#4}) \else
\ifx\six\empty {\rm #1}_{#5}(#3_{#4}) \else
\ifx\seven\empty {\rm #1}_{#5} ({\goth#2})\else                                               
\ifx\eight\empty {\rm #1}_{#5}({#3})                                               
\fi \fi \fi \fi \fi \fi \fi \fi \fi}
\def\hec#1#2#3#4#5{\def\un{#3#4#5}\def\deux{#3#5}\def\trois{#3}
\def\four{#2#4#5}\def\five{#2#5}\def\six{#2}\def\seven{#3#4}
\def\eight{#2#4} \def\nine{#2#3#4}
\ifx\nine\empty \hat{{\rm #1}}_{#5} \else
\ifx\un\empty \hat{{\rm #1}}({\goth #2}) \else
\ifx\deux\empty \hat{{\rm #1}}({\goth #2}_{#4}) \else
\ifx\trois\empty \hat{{\rm #1}}_{#5}({\goth #2}_{#4}) \else
\ifx\four\empty \hat{{\rm #1}}(#3) \else
\ifx\five\empty \hat{{\rm #1}}(#3_{#4}) \else
\ifx\six\empty \hat{{\rm #1}}_{#5}(#3_{#4}) \else 
\ifx\seven\empty \hat{{\rm #1}}_{#5} ({\goth#2})  \else                                             
\ifx\eight\empty \hat{{\rm #1}}_{#5}({#3}) 
\fi \fi \fi \fi \fi \fi \fi \fi \fi}
\def\e#1#2{\ec {#1}#2{}{}{}}
\def\es#1#2{\ec {#1}{}{#2}{}{}}

